\documentclass{amsart}
\usepackage{amsmath,amsthm,amssymb}
\usepackage[latin1]{inputenc}
\usepackage[dvips]{graphicx}

\newtheorem{theo}{\em Theorem}
\newtheorem{prop}{\em Proposition}

\newtheorem{rem}{\em Remark}
\newtheorem{defi}{\em Definition}
\newtheorem{lem}{\em Lemma}
\newtheorem{ex}{\em Example}
\newtheorem{conj}{\em Conjecture}

\def\G{\Gamma}
\def\ov{\overline}
\def\tr{\mathrm{tr}}
\def\cy{\mathrm{cy}}
\def\d{\mathrm{d}}

\def\div{\mathrm{div}}

\newcommand{\R}{\mathbb{R}}
\newcommand{\C}{\mathbb{C}}

\newcommand{\fin}{\hfill $\blacksquare$\\}

\newcommand{\lie}{\mathfrak{lie}}

\newcommand{\krv}{\mathfrak{krv}}
\newcommand{\tder}{\mathfrak{tder}}
\newcommand{\sder}{\mathfrak{sder}}

\newcommand{\der}{\mathfrak{der}}
\newcommand{\ch}{\mathrm{ch}}
\newcommand{\duf}{\mathrm{duf}}
\renewcommand{\t}{\mathfrak{t}}

\begin{document}

\title[]{Kontsevich deformation quantization and flat connections}

\author{Anton Alekseev}
\address{Section de math\'ematiques, Universit\'e de Gen\`eve, 2-4 rue du Li\`evre,
c.p. 64, 1211 Gen\`eve 4, Switzerland}
\email{alekseev@math.unige.ch}

\author{Charles Torossian}
\address{Institut Math\'ematiques de Jussieu, Universit\'e Paris 7, CNRS;
Case 7012, 2 place Jussieu, 75005 Paris, France}
\email{torossian@math.jussieu.fr}

\bibliographystyle{alpha}

\maketitle

\begin{abstract}
In \cite{To1}, the second author used the Kontsevich deformation quantization
technique to define a natural connection $\omega_n$ on the compactified
configuration spaces $\ov{C}_{n,0}$ of $n$  points on the upper half-plane.
Connections $\omega_n$ take values in the Lie algebra of derivations of the
free Lie algebra with $n$ generators. In this paper, we show that $\omega_n$
is flat.

The configuration space $\ov{C}_{n,0}$ contains a boundary stratum at infinity
which coincides with the (compactified) configuration space of $n$ points on
the complex plane. When restricted to this stratum, $\omega_n$ gives rise
to a flat connection $\omega_n^\infty$. We show that the parallel transport
$\Phi$ defined by the connection $\omega_3^\infty$ between configuration
$1(23)$ and $(12)3$ verifies axioms of an associator.

We conjecture that $\omega_n^\infty$ takes values in the Lie algebra $\t_n$
of infinitesimal braids. If correct, this conjecture implies that
$\Phi \in \exp(\t_3)$ is a Drinfeld's associator.
Furthermore, we prove  $\Phi \neq \Phi_{KZ}$ showing that $\Phi$ is a new explicit
solution of the associator axioms.

\end{abstract}
\section{Introduction}

The Kontsevich proof of the formality conjecture and the construction of the
star product on $\mathbb{R}^d$ equipped with a given Poisson structure
make use of integrals of certain differential forms over compactified
configuration spaces $\ov{C}_{n,m}$ of points on the upper half-plane. Here
$n$ points are free to move in the upper half-plane,  $m$ points are bound
to the real axis, and we quotient by the diagonal action of the
group $z \mapsto az +b$ with $a \in \R_+, b\in \R$.

In this paper, we use the same ingredients to study a certain connection $\omega_n$
on $\ov{C}_{n,0}$ with values in the Lie algebra of derivations of the
free Lie algebra with $n$ generators. This connection
was introduced by the second author in \cite{To1}. One of our
results is flatness of $\omega_n$.

The compactified configuration space $\ov{C}_{n,0}$ contains a boundary stratum
``at infinity'' which coincides with the configuration space of $n$ points
on the complex plane (quotient by the diagonal action
$z \mapsto az +b$ with $a \in \R_+, b\in \C$).
Over this boundary stratum, the connection $\omega_n$
restricts to the connection $\omega_n^\infty$ with values in the Lie algebra
$\krv_n$ defined in \cite{AT1}. We conjecture that in fact $\omega_n$
takes values in the Lie algebra $\t_n \subset \krv_n$ defined by the infinitesimal
braid relations.

Let $\Phi$ be the parallel transport defined by the connection
$\omega_3^\infty$
for the straight path between configurations $1(23)$ and $(12)3$ of
3 points on the complex plane. We show that $\Phi$ verifies axioms
of an associator with values in the group $KRV_3=\exp(\krv_3)$. If the conjecture
of the previous paragraph holds true, then $\Phi \in \exp(\t_3)$, and it becomes
a Drinfeld's associator.  The key ingredient in the proof of the
pentagon axiom if the flatness property of $\omega_n^\infty$.
The construction of $\Phi$ is parallel to the construction of the
Knizhnik-Zamolodchikov associator $\Phi_{KZ}$ in \cite{Dr} with $\omega_3^\infty$
replacing the Knizhnik-Zamolodchikov connection. Furthermore, one can show that
$\Phi$ is even, and hence $\Phi \neq \Phi_{KZ}$.

While this paper was in preparation, we learnt of the work \cite{SW}
proving our conjecture stated above.

The plan of the paper is as follows. In Section 2, we review some standard facts
about the Kontsevich déformation quantization technique and  free Lie algebras.
In Section 3, we prove flatness of the connection $\omega_n$. Section 4 contains
the proof of associator axioms for the element $\Phi$.\\

{\bf Acknowledgements:}
We thank D. Barlet and F. Brown for useful discussions and remarks.
We are grateful to P. Severa and T. Willwacher for letting us know of their
forthcoming work \cite{SW}.
Research of A.A. was supported in part by the grants
200020-121675 and 200020-120042 of the Swiss National Science Foundation.
Research of C.T. was supported by CNRS.

\section{Deformation quantization and free Lie algebras}

Many sources are now available on the Kontsevich formula for quantization
of Poisson brackets (see {\em e.g.}  \cite{CKT}). For convenience of the reader,
we briefly recall the main ingredients of \cite{Kont} for $\R^d$  and
the construction \cite{To1} of the connection $\omega_n$.

\subsection{Free Lie algebras and their derivations}

\subsubsection{Free Lie algebras and derivations}
Let $\mathbb{K}$ be a field of characteristic zero, and let $\lie_n=\lie(x_1,\dots, x_n)$
be the degree completion of the graded free Lie algebra over $\mathbb{K}$  withgenerators
$x_1, \dots, x_n$ of degree one. We shall denote by $\der_n$ the Lie algebra of derivations
of $\lie_n$. An element $u \in \der_n$ is completely determined by its
values on generators, $u(x_1), \dots, u(x_n) \in \lie_n$.
The Lie algebra $\der_n$ carries a grading induced by the one of $\lie_n$.

\begin{defi}
A derivation $u \in \der_n$ is called {\em tangential} if there exist
$a_i \in \lie_n, i=1, \dots, n$ such that $u(x_i)=[x_i, a_i]$.
\end{defi}
Tangential derivations form a Lie subalgebra $\tder_n \subset \der_n$.
Elements of $\tder_n$ are in one-to-one correspondence with $n$-tuples
of elements of $\lie_n$, $(a_1, \dots, a_n)$, which verify the condition
that $a_k$ has no linear term in $x_k$ for all $k$. By abuse of notations,
we shall often write $u=(a_1, \dots, a_n)$.  For two elements of $\tder_n$,
$u=(a_1, \dots, a_n)$ and $v=(b_1, \dots, b_n)$, we have
$[u,v]_{\tder}=(c_1,\dots,c_n)$ with
\begin{equation} \label{3terms}
c_k=u(b_k) - u(a_k) + [a_k, b_k]_{\lie} .
\end{equation}

\begin{defi}
A derivation $u = (a_1, \dots, a_n)\in \tder_n$ is called {\em special} if
$u(x)=\sum_i [x_i, a_i]=0$ for $x=\sum_{i=1}^n x_i$.
\end{defi}
We shall denote the space of special derivations by $\sder_n$.
It is obvious that $\sder_n \subset \tder_n$ is a Lie subalgebra.
Both $\tder_n$ and $\sder_n$ integrate to prouniportent groups
denoted by $TAut_n$ and $SAut_n$, respectively. In more detail,
$TAut_n$ consists of automorphisms of $\lie_n$ such that
$x_i \mapsto {\rm Ad}_{g_i} x_i = g_i x_i g_i^{-1}$, where $g_i \in \exp(\lie_n)$.
Similarly, elements of $SAut_n$ are tangential automorphisms of $\lie_n$
with an extra property  $x=\sum_{i=1}^n x_i \mapsto x$.

The family of Lie algebras $\tder_n$ is equipped with simplicial
Lie homomorphisms $\tder_n \rightarrow \tder_{n+1}$. For instance,
for $u=(a,b) \in \tder_2$ we define
$$
\begin{array}{lll}
u^{1,2} & = & (a(x,y), b(x,y), 0) \, , \\
u^{2,3} & = & (0, a(y,z), b(y,z)) \, , \\
u^{12,3} & = & (a(x+y, z), a(x+y,z), b(x+y,z)) ,
\end{array}
$$
and similarly for other simplicial maps.
These Lie homomorphisms integrate to group homomorphisms of $TAut_n$ and
$SAut_n$.

\subsubsection{Cyclic words}  \label{subsec:tr}
Let $Ass^+_n=\prod_{k=1}^\infty Ass^k(x_1, \dots, x_n)$ be the graded
free associative algebra (without unit) with generators $x_1, \dots, x_n$.
Every element $a \in Ass^+_n$
admits a unique decomposition of the form $a=\sum_{i=1}^n (\partial_i a) x_i$,
where $a_i \in Ass_n$ ($Ass_n$ is a free associative algebra with unit).

We define the graded vector space $\cy_n$ as a quotient
$$
\cy_n = Ass^+_n/\langle (ab-ba); a,b \in Ass_n \rangle .
$$
Here $\langle (ab-ba); a,b \in Ass_n \rangle$
is the subspace of $Ass^+_n$ spanned by commutators. The multiplication map of $Ass^+_n$
does not descend to $\cy_n$ which only has a structure of a graded vector
space. We shall denote by $\tr: Ass^+_n \to \cy_n$ the natural
projection. By definition, we have $\tr(ab)=\tr(ba)$ for all $a,b \in Ass_n$
imitating the defining property of trace. In general, graded components
of $\cy_n$ are spanned by words of a given length modulo cyclic permutations.

\begin{ex}
The space $\cy_1$ is isomorphic to the space of formal power series
in one variable without constant term, $\cy_1 \cong x k[[x]]$.
This isomorphism is given by the following formula,
$$
f(x)=\sum_{k=1}^\infty f_k x^k \mapsto  \sum_{k=1}^\infty f_k \tr(x^k) .
$$
\end{ex}
\subsubsection{Divergence}
Let $ u=(a_1, \dots, a_n) \in \tder_n$. We define the divergence as
$$
\div (u)= \sum_{i=1}^n \tr(x_i (\partial_i a_i)).
$$
It is a 1-cocycle of $\tder_n$ with values in $\cy_n$
(see Proposition 3.6 in \cite{AT1}).

We define $\krv_n\subset \sder_n \subset \tder_n $ as
the Lie algebra of special derivation with vanishing divergence.
Hence, $u=(a_1,\ldots, a_n) \in \krv_n$ is a solution of
two equations:  $\sum_{i=1}^n [x_i, a_i]=0$ and
$\sum_{i=1}^n \tr(x_i (\partial_i a_i))=0$.
We shall denote by $KRV_n=\exp(\krv_n)$ the
corresponding prounipotent group.

\subsection{Kontsevich construction}
\subsubsection{Configurations spaces}
We denote by $C_{n,m}$ the configuration space of  $n$ distincts points in
the upper half plane and $m$ points on the real line modulo the diagonal
action of the group  $z \mapsto az+b$ ($a \in \R_+, b\in \R$).
In \cite{Kont},
Kontsevich constructed compactifications of spaces $C_{n,m}$ denoted by
$\overline{C}_{n,m}$. These are manifolds with corners
of dimension  $2n-2+m$. We denote by  $\overline{C}^{+}_{n,m}$ the
connected component of $\overline{C}_{n,m}$
with real points in the standard order (\textit{id.}
$\ov{1}< \ov{2}<\cdots< \ov{m}$).

\begin{figure}[h!]
\begin{center}
\includegraphics[width=5cm]{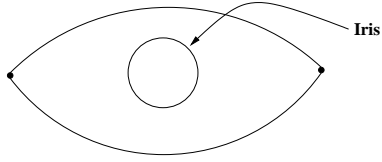}
\caption{\footnotesize Variety  $\overline{C}_{2, 0}$.}\label{oeil.eps}
\end{center}
\end{figure}

The compactified configuration space $\overline{C}_{2, 0}$ (the ``Kontsevich eye'')
is shown on Fig.~\ref{oeil.eps}. The upper and lower eyelids correspond to one of the
points ($z_1$ or $z_2$) on the real line, left and right corners of the eye
are configurations with $z_1,z_2 \in \R$ and $z_1>z_2$ or $z_1<z_2$. The boundary
of the iris takes into account configurations where $z_1$ and $z_2$ collapse
inside the complex plane. The angle along the iris keeps track of the
angle at which $z_1$ approaches $z_2$.

\subsubsection{Graphs}

A graph $\Gamma$ is a collection of vertices $V_\Gamma$ and oriented
edges $E_\G$. Vertices are ordered, and the edges are ordered in a way compatible
with the order of the vertices.  We denote by $G_{n,2}$ the set of graphs with $n+2$
vertices and $2n$ edges verifying the following properties:\\

\smallskip
i - There are $n$ vertices of the first type $1,2,
\cdots,  n$ and $2$ vertices of the second type $\overline{1}, \overline{2}$

ii - Edges start from vertices of the first type, 2 edges per vertex.

iii - Source and target of an edge are distinct.

iv - There are no multiple edges (same source and target). \\

\begin{figure}[!h]
\begin{center}
\includegraphics[width=7cm]{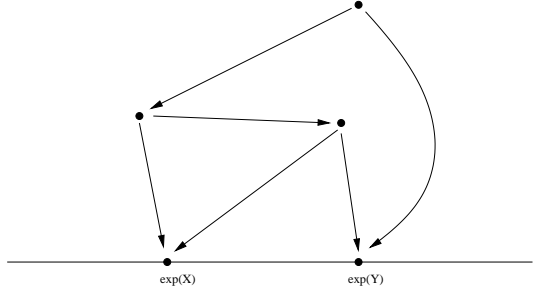}
\caption{\footnotesize Lie type graph with symbol
$\G(x, y)=[[x, [x, y]], y]$.}\label{Lie.eps}
\end{center}
\end{figure}

\begin{figure}[h!]
\begin{center}
\includegraphics[width=6cm]{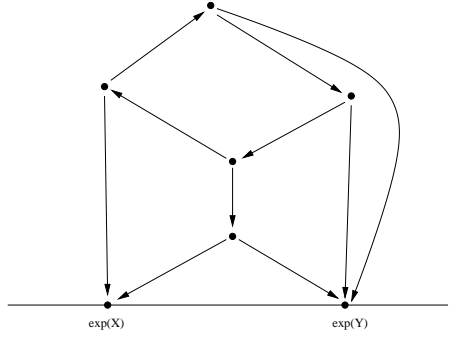}\caption{\footnotesize
Wheel type graph with symbol $ \Gamma(x,y)=\tr(y^2[x,y]x).$}
\label{roue.eps}
\end{center}
\end{figure}

We are interested in the case of \textit{linear graphs}. That is,
vertices of the first type admit at most one incoming edge. Such
graphs are superpositions of Lie type graphs (graphs with one root as
on Fig.~\ref{Lie.eps}) and  wheel type graphs (graph with one
oriented loop, as on Fig.~\ref{roue.eps}).

\subsubsection{The angle map and  Kontsevich weights}

Let $p$ and $q$ be two points on the upper half plane. Consider the hyperbolic
angle map on $C_{2,0}$:
\begin{equation}
\phi_h(p,q)=\arg\left( \frac{q-p}{q-\overline{p}}\right) \in \mathbb{T}^1.
\end{equation}
This function admits a continuous extension to the compactification
$\overline{C}_{2,0}$.

Consider a graph $\Gamma \in G_{n,2}$, and draw it in the upper half plane with
vertices of the second type on the real line. By restriction, each edge $e$ defines
an angle map  $\phi_{e}$ on  $\overline{C}^+_{n,2}$. The ordered product
\begin{equation}
\Omega_{\Gamma}=\bigwedge _{e \in E_{\Gamma}}\d\phi_{e}
\end{equation}
is a regular  $2n$-form on $\overline{C}^+_{n,2}$ (which is a  $2n$-dim compact space).
\begin{defi} The Kontsevich weight of  $\G$ is given by the
following formula,
\begin{equation}
w_{\Gamma}=\frac{1}{(2\pi)^{2n}}\int_{\overline{C}^+_{n,2}} \Omega_{\Gamma}.
\end{equation}
\end{defi}

\subsection{Campbell-Hausdorff and Duflo formulas}\label{sectionKVnew}
Lie type graphs in $G_{n,2}$ are binary rooted trees. Hence, to
each $\Gamma \in G_{n,2}$ one can associate a Lie word
$\Gamma(x, y)\in \lie_2$ of degree $2n$ in variables $x, y$
(see Fig.~\ref{Lie.eps}).
Similarly, if $\G$ is a  wheel type graph, it corresponds to an
element $G(x,y) \in \cy_2$  (see Fig.~\ref{roue.eps}).

Recall the definition of the Duflo density function
$$
\duf(x,y)=\frac{1}{2} \left(j(x)+j(y) - j(\ch(x,y))\right) \in \cy_2 ,
$$
where $\ch(x,y)=\log(e^xe^y)$ is the Campbell-Hausdorff series and
$$
j(x)=\sum_{n\geq 2} \frac{b_n}{n \cdot n!} \tr(x^n)
$$
with $b_n$ the Bernoulli numbers. The following Theorem
relates functions $\ch(x,y)$ and $\duf(x,y)$ to the Kontsevich
graphical calculus.

\begin{theo}[\cite{Ka}, \cite{AST}]\label{theoKVnew}
The following identities hold true:
\begin{equation} \ch(x,y)=x+y +\sum\limits_{n\geq
1}\sum\limits_{\substack{
 \Gamma  \;\mathrm{simple}\\\mathrm{geometric}\\
\mathrm{Lie \;type}\; (n,2)}}w_{\Gamma} \Gamma(x,y),
\end{equation}
\begin{equation}
\duf(x,y)=\sum\limits_{n\geq 1}\sum\limits_{\substack{
\Gamma \;\mathrm{simple}\\\mathrm{geometric}\\
\mathrm{wheel\; type\; \; (n,2)}}} \frac{ w_{\Gamma}
}{m_\Gamma}\G(x,y),
\end{equation}
where $m_\G$ is the order of the symmetry group of the graph $\G$.
\end{theo}
\noindent Here \textit{geometric} means that graphs are not labeled. Note that
the definition of both $\G(x,y)$ and $w_\G$ requires an order on the set
of edges, but the product $w_\G \G(x,y)$ is independent of this order.
Even though $\ch(x,y)$ and $\duf(x,y)$ are defined over rationals, some of the
coefficients $w_\G$ are very probably irrational (see example of~\cite{FW}).

\subsection{$\xi$-deformation}\label{sectiondeformationkont}
In  \cite{To1},  one studies the following deformation for
the Campbell-Hausdorff formula. Let $\xi \in \overline{C}_{2, 0}$,
$\G\in G_{n,2}$, and let $\pi$ be the natural projection
from $\overline{C}_{n+2, 0}$  onto $\ov{C}_{2,0}$. We define
the coefficients $w_\Gamma(\xi)$ for $\xi \in \ov{C}_{2, 0}$  as
$$
w_\Gamma(\xi)=\frac{1}{(2\pi)^{2n}}\int_{\pi^{-1}(\xi)} \Omega_\Gamma .
$$
Functions $w_\Gamma(\xi)$ are smooth over $C_{2,0}$, and they are  continuous
over the compactification $\ov{C}_{2, 0}$.
The $\xi$-deformation of the Campbell-Hausdorff series $\ch_\xi(x, y)$  is
defined as
\begin{equation}\label{BCH}\ch_\xi(x, y)= x+y +\sum\limits_{n\geq
1}\sum\limits_{\substack{
 \Gamma  \;\mathrm{simple}\\\mathrm{geometric}\\
\mathrm{Lie \;type}\; (n,2)}} w_\G(\xi) \G(x, y).
\end{equation}
In a similar fashion, we introduce a deformation of the Duflo function,
$$
\duf_\xi(x,y)=\sum\limits_{n\geq 1}\sum\limits_{\substack{
\Gamma \;\mathrm{simple}\\\mathrm{geometric}\\
\mathrm{wheel\; type\; \; (n,2)}}} \frac{ w_{\Gamma}(\xi)}{m_\Gamma}\G(x,y) .
$$
For $\xi=(0,1)$ (the right corner of the eye on Fig.~\ref{oeil.eps}),
the expression (\ref{BCH}) is  given by the standard
Campbell-Hausdorff series, and for $\xi$ in the position $\alpha$
on the iris, the Kontsevich Vanishing Lemma implies $\ch_\alpha(x, y)=x+y$.
By the results of \cite{Shoi}, for $\xi=(0,1)$ the Duflo function
$\duf_\xi(x,y)$ coincides with the standard Duflo function,
and for $\xi$ on the iris one has $\duf_\alpha(x,y)=0$.

\subsection{Connection $\omega_2$}\label{defiOmega2}

In \cite{To1}, one defines a  connection on ${C}_{2,0}$ with values
in $\tder_2$,
$$
\omega_2=\big(F_\xi(x,y), G_\xi(x,y)\big) .
$$
Here $F_\xi$ and $G_\xi$ are $1$-forms on ${C}_{2,0}$ taking values in $\lie_2$.
They satisfy the following two (Kashiwara-Vergne type) equations
(see Theorem 1 and Theorem 2 in \cite{To1})
\begin{equation}\label{KVcharles1}
\d \, \ch_\xi (x, y)= \omega_2(\ch_\xi(x,y))
\end{equation}
\begin{equation}\label{KVcharles2}
\d \, \duf_{\xi}(x, y) = \omega_2(\duf_\xi(x,y)) + \div(\omega_2),
\end{equation}
where $\d$ is the de Rham differential on $C_{2,0}$, and
$\omega_2$ acts on $\ch_\xi(x,y)$ and $\duf_\xi(x,y)$ as
a derivation of $\lie_2$.

\begin{figure}[h!]
\begin{center}
\includegraphics[width=5cm]{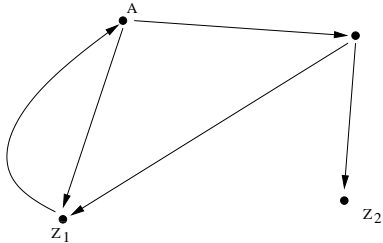}
\caption{\footnotesize Extended graph}\label{extended.eps}
\end{center}
\end{figure}

Let us briefly recall the construction of $\omega_2$.
We will denote by  $A, B\in G_{n,2}$ connected graphs of Lie type,
and we define an extended graph $ \Rsh \!\! A$ (resp.  $B\Lsh$ ) as a graph
with an additional edge starting at $\ov{1}$ (resp. $\ov{2}$) and ending at
the root of $A$ (resp. $B$), see Fig.~\ref{extended.eps}.
Note that $n=0$ is allowed, but since the source and the target of an edge must be
distinct, $\Gamma$ will have a single edge starting at $\ov{1}$ and ending at
$\ov{2}$, or starting at $\ov{2}$ and ending at $\ov{1}$.

Draw the extended graph in the upper half plane
(with vertices of the second type corresponding to $\xi \in \ov{C}_{2,0}$). Then,
$$
\Omega_{\Rsh A}= \bigwedge _{e \in E_{\Rsh\!\! A}}\d\phi_{e}
$$
is a $2n+1$-form  on $\ov{C}_{n+2, 0}$ (which is a $2n+2$-dim compact space).
The push forward along the natural projection
$\pi : \ov{C}_{n+2, 0} \rightarrow \ov{C}_{2, 0}$ yields a 1-form
on $\ov{C}_{2, 0}$, $\omega_{\Rsh A}= \pi_*(\Omega_{\Rsh A})$.

The connection $\omega_2= (F_\xi(x, y), G_\xi(x, y))$ is
defined by the following formula\footnote{For $n=0$, $A=y$ and $B=x$.},
\begin{equation}
\left\{\begin{array}{ccc}
  F_\xi(x, y) & = & \sum\limits_{n\geq 0} \sum\limits_{\substack{A\,\mathrm{ graph \, of\,}\\\mathrm{ Lie\, type} \; (n,2)}} \omega_{\Rsh
A} A(x, y) \\
  G_\xi(x,y) &= &\sum\limits_{n\geq 0} \sum\limits_{\substack{B, \mathrm{graph \, of\,}\\\mathrm{ Lie\, type} \; (n,2) }}\omega_{ B\Lsh} B(x, y)
\end{array}\right.
\end{equation}

\subsection{Definition of $\omega_n$}\label{defiOmegan} We now extend
this construction  to an arbitrary number of vertices of the second type.

Consider $\Gamma\in G_{p, n}$ a connected graph of Lie type with vertices
of the second type labeled $\ov{1}, \ldots, \ov{n}$. Define the extended
graph $\Gamma^{(i)}$ by adding an edge from the vertex $\ov{i}$ to the
root of $\G$. Consider the natural projection $\pi: \ov{C_{p+n, 0}}
\rightarrow \ov{C_{n, 0}}$ and
take the push forward $1$-form $\omega_{\G^{(i)}}=\pi_*(\Omega_{\G^{(i)}})$.
We define 1-forms with values in $\lie_n$
$$
F_i=  \sum\limits_{p\geq 0} \sum\limits_{\substack{\Gamma\,\mathrm{ graph \,
of\,}\\\mathrm{ Lie\, type} \; (p,n)}}\omega_{{\G}^{(i)}} \Gamma(x_1,\ldots,  x_n),
$$
and  $\omega_n=(F_1, \ldots, , F_n)$ yields a connection with values
in $\tder_n$.

The connection $\omega_n$ is smooth over $C_{n,0}$. Over the compactification
$\ov{C}_{n,0}$, it belongs  to the class $L^1$ when restricted
to piece-wise differentiable curves. Hence, along such curves
all iterated integrals converge, and there is a unique solution of the
initial value problem $\d g = -g \omega$ with
$g(z_0)=1$ for the base point $z_0$ (e.g. by using the
Gr\"onwall's Lemma). Therefore, parallel transports
are well defined.

The same applies to restrictions of
$\omega_n$ to boundary strata of $\ov{C}_{n,0}$ of dimension at least one.
For instance, in the case of $\ov{C}_{2,0}$ one can consider a path along the eyelid,
or a generic path from the corner of $\ov{C}_{2, 0}$ to the iris.

We will need restrictions of $\omega_n$ to various boundary strata of
co-dimension one of $\ov{C}_{n,0}$. First of all, there is a stratum ``at infinity''
equal to the configuration space $C_n$ of $n$ points on the complex plane
(modulo the diagonal action of the group $z \mapsto az+b$ for $a\in \R_+, b \in \C$).
We denote the corresponding connection $\omega_n^\infty$. It is given by the same
formula as $\omega_n$ with the configuration space $\ov{C}_{n,0}$ replaced
by $\ov{C}_n$, and the hyperbolic angle replaced with the Euclidean angle.

Next, for $q$ points collapsing inside the upper half plane, we have a stratum
of the form $C_q \times C_{n-q+1,0}$. We denote the natural projections by $\pi_1$
and $\pi_2$, and obtain an expression for the connection
$$
\omega_n|_{C_q \times C_{n-q+1,0}}= \pi_1^* (\omega^\infty_q)^{1,2,\dots,q} +
\pi_2^* \omega_{n-q+1}^{12\dots q,q+1,\dots,n} .
$$
A similar property holds for the connection $\omega_n^\infty$ for the stratum
$C_q \times C_{n-q+1}$ corresponding to $q$ points collapsing together,
$$
\omega^\infty_n|_{C_q \times C_{n-q+1}}= \pi_1^* (\omega^\infty_q)^{1,2,\dots,q} +
\pi_2^* (\omega_{n-q+1}^\infty)^{12\dots q,q+1,\dots,n} .
$$

In the case when $q$ points are collapsing to the point on the real axis,
we obtain the stratum $C_{q,0} \times C_{n-q,1}$, and for the connection
we get
$$
\omega_n|_{C_{q,0} \times C_{n-q,1}} = \pi_1^* \omega_q^{1,2,\dots,q}
+\pi_2^* \omega_{n-q+1}^{12\dots q,q+1,\dots,n}|_{C_{n-q,1}} .
$$
Note that the restriction of the connection form $\omega_n$ to the
boundary stratum $C_{n-1,1}$ corresponds to configurations with the point $z_1$
on the real axis, and it has
the following property: its first component (as an element of $\tder_n$)
vanishes since the 1-form $d\phi_e$ vanishes when the source of the edge $e$
is bound to the real axis.


\section{Zero curvature equation and applications}\label{zero}

One of our main results is flatness of the connection $\omega_n$.

\subsection{The zero curvature equation}

\begin{theo}\label{theocourbure} The connection $\omega_n$ is flat.
That is, the following $2$-form on $C_{n,0}$ vanishes,
\begin{equation}
\d \omega_n + \frac12[\omega_n,\omega_n]=0
\end{equation}
 \end{theo}

\noindent \textit{Proof : } The argument is based on the Stokes formula,
and we give details in the case of $\omega_2$. The case of arbitrary $n$
is treated in a similar fashion.

Let $C_\xi$ be a small circle around $\xi\in C_{2, 0}$,
$\Delta_\xi$ be the corresponding disk, and consider
 $\pi^{-1}(\Delta_\xi)$. Since the forms
$\Omega_{\Rsh A}, \Omega_{B\Lsh}$  are closed, we have
$$
\int_{\pi^{-1}(\Delta_\xi)} \d \,
\Big(\sum\limits_{A} \Omega_{\Rsh A} A(x, y),
\sum\limits_{B}\Omega_{B\Lsh} B(x, y)\Big)=0.
$$
By applying Stokes formula and the definition of the connection $\omega_2$, one obtains
\begin{equation}\label{equationcurvature=0}
0= \int_{C_\xi} \omega_2 + \int\limits_{\bigcup\limits_{z\in \Delta_\xi}
\partial(\pi^{-1}(z))} \Big(\sum\limits_{A} \Omega_{\Rsh A} A(x, y),
\sum\limits_{B}\Omega_{B\Lsh} B(x, y)\Big).
\end{equation}
By using again the Stokes's formula (on the disk $\Delta_\xi$),
one can rewrite the first term
in the form $\int_{\Delta_\xi} \d_{\xi} \omega_2$.

\begin{figure}[h!]
\begin{center}
\includegraphics[width=6cm]{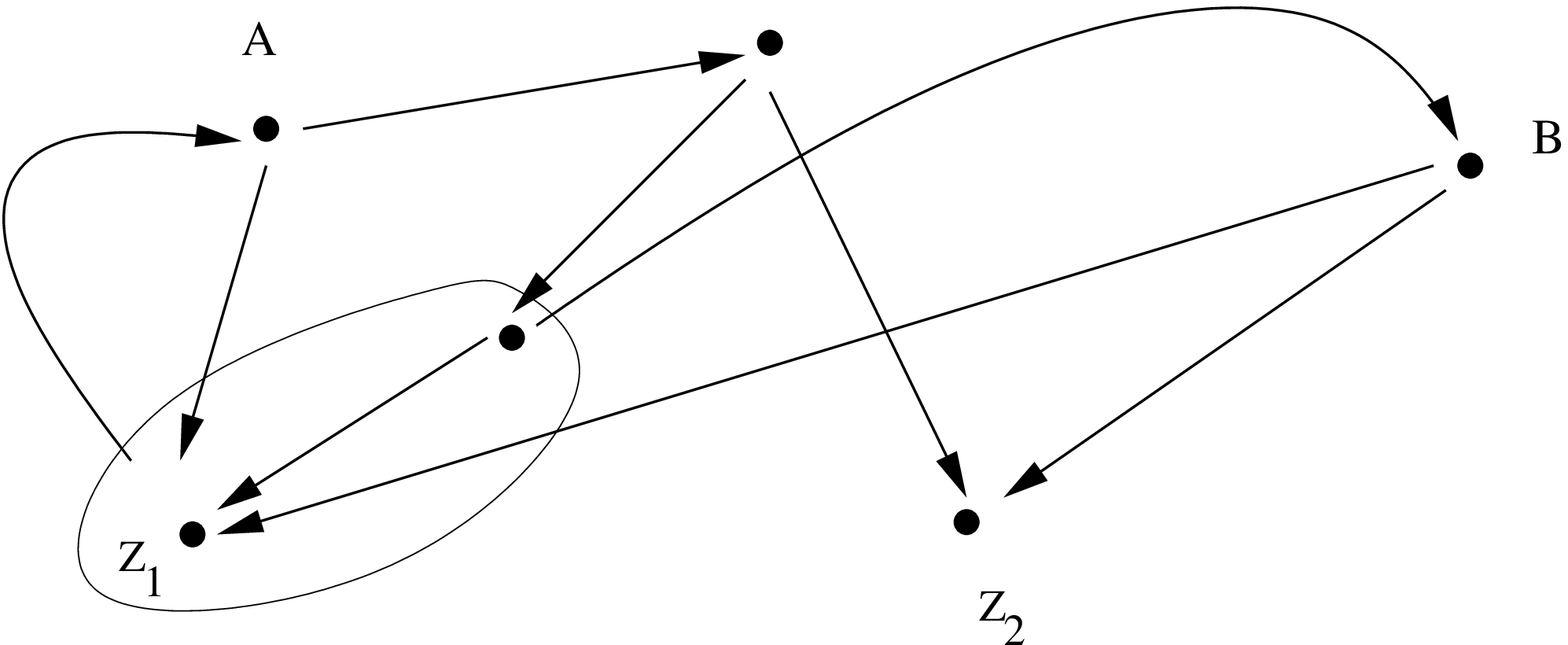}
\caption{\footnotesize}\label{curvature1.eps}
\end{center}
\end{figure}

For the second term, one obtains contributions from the boundary strata of co-dimension one.
The usual arguments in Kontsevich theory rule out strata where more than two points collapse
(by the Kontevich Vanishing Lemma), and strata corresponding to collapse of internal
edges (by Jacobi identity). The remaining strata correspond to collapsing of a vertex
of the first type and a vertex of the second type. Figures below illustrate different cases
of such boundary strata (for the first component):

- Fig.~\ref{curvature1.eps} computes terms of the type
$[x,\omega_{\Rsh B} B(x, y)]\cdot \partial_x\Big(\omega_{\Rsh A}
A(x, y)\Big).$

\begin{figure}[h!]
\begin{center}
\includegraphics[width=6cm]{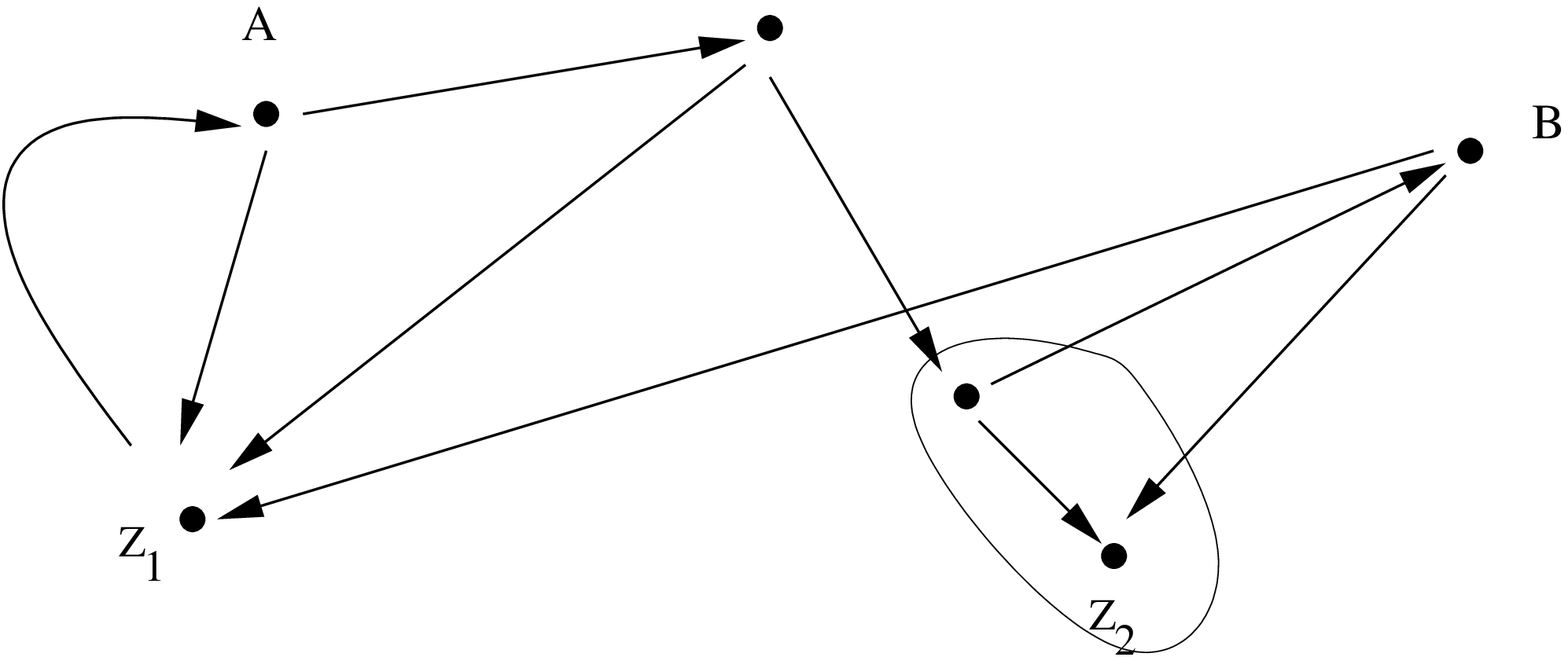}
\caption{\footnotesize}\label{curvature2.eps}
\end{center}
\end{figure}

- Fig.~\ref{curvature2.eps} represents terms of the type $[y,\omega_{
B\Lsh} B(x, y)]\cdot \partial_y \Big(\omega_{\Rsh A} A(x, y)\Big).$

\begin{figure}[h!]
\begin{center}
\includegraphics[width=6cm]{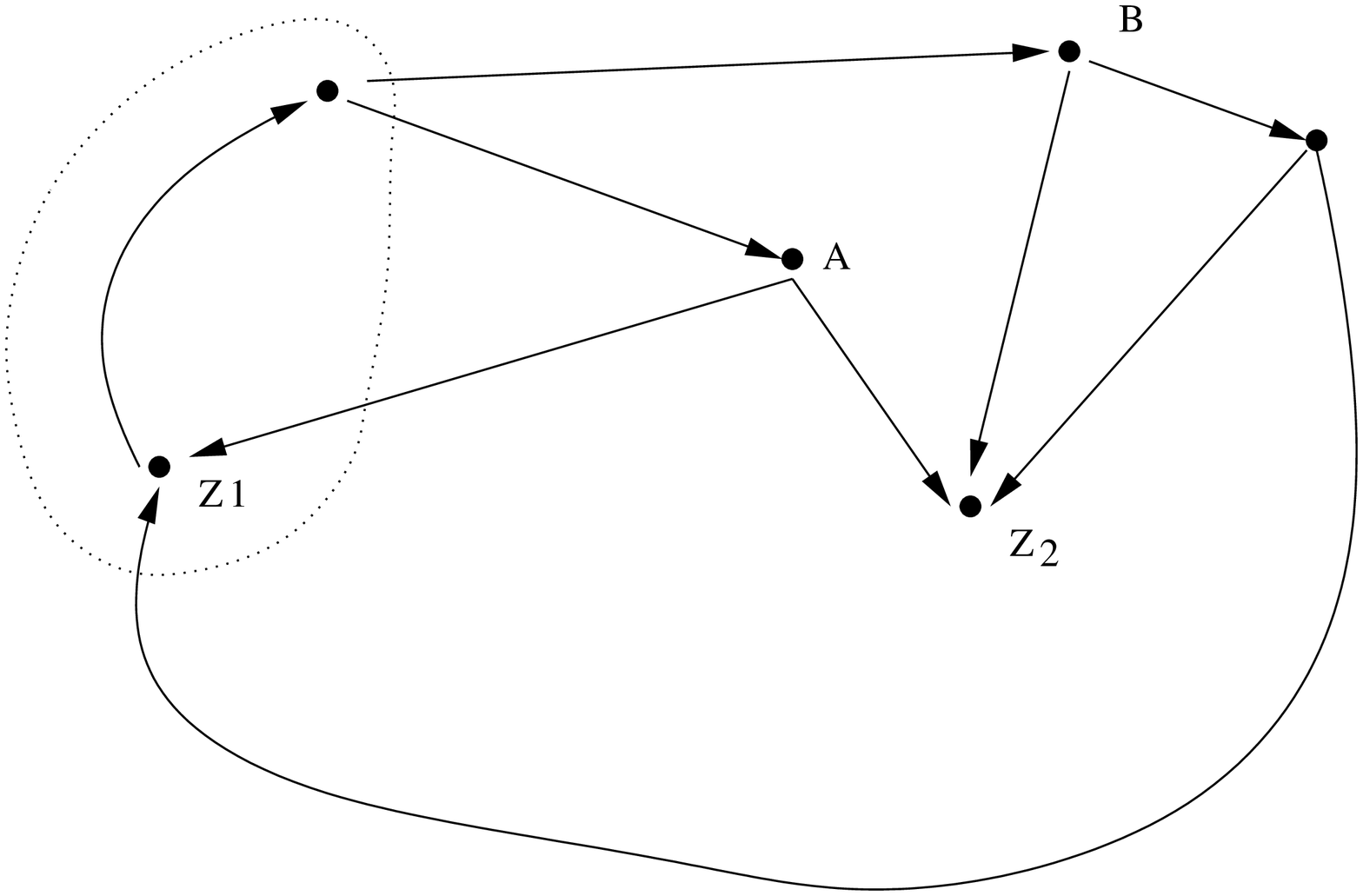}
\caption{\footnotesize} \label{curvature3.eps}
\end{center}
\end{figure}

- Fig.~\ref{curvature3.eps} computes terms of the type $[\omega_{ \Rsh A}
A(x, y), \omega_{ \Rsh B} B(x, y)].$ This term appears only once. It
corresponds to the componentwise bracket in
$\lie_2\times \lie_2$$$\frac12[\omega_2(x, y),
\omega_2(x, y)]_{\lie_2}.$$

These are exactly the three terms of the bracket
$\frac12[\omega_2,\omega_2]_{\tder}$ (see equation (\ref{3terms})).
By equation (\ref{equationcurvature=0}), one gets
$$
\int_{\Delta_\xi} \left(\d \omega_2 + \frac12[\omega_2, \omega_2]\right)=0.
$$
Since the curvature $\d \omega_2 + \frac12[\omega_2, \omega_2]$ is a continuous
function of $\xi$, we conclude that it vanishes on $C_{2,0}$.
\fin

\subsection{Parallel transport and symmetries}\label{sectionKV}

In this Section we discuss various properties of the connection
$\omega_2$, including the induced holonomies and their symmetries.

\subsubsection{Parallel transport}

Since $\omega_2$ is flat, the equation
$$
\d \, g= - g \omega_2,
$$
has a local solution on $C_{2,0}$ with values in
$TAut_2= \exp(\tder_2)$. By abuse of notations we write
$(u,v) \in TAut_2$ for an element acting on generators
by $x\mapsto {\rm Ad}_u x =uxu^{-1}, y\mapsto {\rm Ad}_v y=vyv^{-1}$.

\begin{figure}[h!]
\begin{center}
\includegraphics[width=4cm]{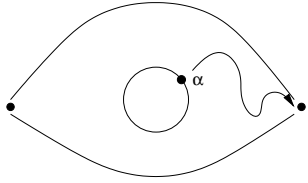}
\caption{\footnotesize A simple path from Duflo-Kontsevich star
product to
 standard product}\label{pathIrisCorner.eps}
\end{center}
\end{figure}

Take  the initial data $g_\alpha=1$ for $\alpha$ on the iris of $\ov{C}_{2, 0}$
(see Fig.~ \ref{pathIrisCorner.eps}), and consider a  path  from
$\alpha$ to $\xi$.  The value at $\xi$ for the parallel transport is
well-defined since the connection is integrable, and by the flatness property
it only depends of the homotopy class of the path. Integrating
the equation (\ref{KVcharles1}), we obtain
$g_\xi(\ch_\xi(x,y))=\ch_\alpha(x,y)=x+y$.

Recall \cite{AM2} that for $\xi=(0,1)$ 
this parallel transport $F$ defines
a solution of the Kashiwara-Vergne conjecture \cite{KV}. We conclude that
this solution is independent of the choice of a path in the trivial homotopy class
(the straight line joining $\alpha=0$ and $\xi=(0,1)$).

\subsubsection{Holonomy}\label{holonomy}

Solutions of equation $\d \, g= - g \omega_2$ are not globally
defined on  $\ov{C}_{2, 0}$ because of the holonomy around the iris.

\begin{lem}\label{restrictionO2} The restriction of $\omega_2$ to the iris is
equal to $\omega_\theta=  \frac{\d \, \theta}{2\pi} (y, x) .$ The  holonomy
around the Iris  $H_{2\pi}$  is given by the inner automorphism
$(\exp(x+y), \exp(x+y)) \in TAut_2$.
\end{lem}

\noindent  \textit{Proof : }  When the point $\xi \in \ov{C}_{2,0}$
reaches the iris, the Kontsevich angle map degenerates to the Euclidian angle
map on the complex plan, and the connection $\omega_2$ is replaced by $\omega_2^\infty$
$$
\omega_2^\infty = \Big(\sum w^\infty_{\Rsh A} A(x, y),\sum w^\infty_{ B\Lsh} B(x, y)\Big).
$$
Since the Euclidian angle  is rotation invariant,
so is the $1$-form $w^\infty_{\Rsh A}$. Hence, it is sufficient
to compute $\int_{\mathbb{T}^1} w^\infty_{\Rsh A}$. By
Kontsevich Vanishing Lemma (see \cite{Kont} \S~6.6), integrals of
$3$ and more angle $1$-forms vanish. Therefore, for $A$ a nontrivial
graph one gets $\int_{\mathbb{T}^1} w^\infty_{\Rsh A} =0$ which
implies
$w^\infty_{\Rsh A} =0$. As a result, we obtain the connection
$\omega_\theta$ by adding
two trivial graph contributions,
$$
\omega_\theta= \frac{\d \, \theta}{2\pi} ( y, 0) +
\frac{\d \, \theta}{2\pi} ( 0, x) = \frac{\d \, \theta}{2\pi} ( y, x).
$$
Let's integrate the equation $\d_\theta g= -g \omega_\theta$ over
the boundary of the Iris.
Note that $t=(y,x)$ is actually an inner derivation since
$t(x)=[x,y]=[x,x+y]$ and $t(y)=[y,x]=[y,x+y]$. We conclude that the parallel
transport around the iris  is given by
$H_\theta=\exp(\theta t/2\pi)=(\exp(\theta(x+y)/2\pi), \exp(\theta(x+y)/2\pi)$.
In particular, for $\theta=2\pi$ we obtain
$H_{2\pi}=(\exp(x+y), \exp(x+y))$, as required.
\fin

\subsubsection{Symmetries of the connection}

\noindent Consider the following involutions on $\ov{C}_{2, 0}$,
$$
\sigma_1 : (z_1, z_2) \mapsto (z_2, z_1) \quad \mathrm{and}
\quad \sigma_2 : (z_1, z_2) \mapsto (-\bar{z}_1, -\bar{z}_2).
$$
Identifying $\ov{C}_{2, 0}$ with the Kontsevich eye (see Fig.~\ref{oeil.eps}),
$\sigma_1$ is the reflection with respect to the center of the eye, and
$\sigma_2$ is the reflection with respect to the vertical axis
(see Fig.~\ref{oeil.eps}).

We shall denote by $\tau_1$ and $\tau_2$ the following involutions of
$\tder_2$,
$$
\begin{array}{lll}
\tau_1 & : & (F(x, y), G(x, y))\mapsto (G(y, x), F(y, x)),  \\
\tau_2 & : & (F(x,y), G(x,y)) \mapsto (F(-x,-y), G(-x,-y)).
\end{array}
$$
They lift to involutions of $TAut_2$.

\begin{prop} \label{symmetries}
 The connection $\omega_2$ verifies
$$
\sigma^*_1(\omega_2)=\tau_1(\omega_2)
\hskip 0.3cm  , \hskip 0.3cm
\sigma^*_2(\omega_2)=\tau_2(\omega_2)
$$
\end{prop}

\noindent \textit{Proof : }
The involution $\sigma_1$ simply exchanges the colors $x$ and $y$ of all
graphs which induces the involution $\tau_1$ on $\tder_2$.

The involution $\sigma_2$ flips the sign of the one form $d \phi_e$ for
each edge (since the reflection
changes sign of the Euclidean angle), and changes the orientation of
each integration over a complex variable. Hence, for a graph  with $n$
internal vertices we collect $-1$ to the power $(2n+1)+n \equiv n+1 \,\, (mod \, \, 2)$.
Corresponding rooted trees have exactly $n+1$ leaves. Hence, one should change a sign
of each leaf which results in applying the  involution $\tau_2$.
\fin

Let $F\in TAut_2$ be the parallel transport of the equation
$\d g=-g \omega_2 $ for the straight path between the
the position $0$ on the iris to the right corner of the eye
$\ov{C}_{2,0}$. Since the path is invariant under the composition
$\sigma_1 \circ \sigma_2=\sigma_2 \circ \sigma_1$,
the parallel transport is invariant under $\tau=\tau_1 \tau_2=\tau_2\tau_1$,
$\tau(F)=F$.

In order to discuss the involution $\tau_1$ and $\tau_2$ separately,
we need the following Lemma.

\begin{lem} The parallel transport  along the lower eyelid in the
counter-clockwise direction is equal to $R=(\exp(y),1) \in TAut_2$.
\end{lem}

\noindent \textit{Proof : }
The connection restricted to the lower eyelid has a trivial second component
because edges starting from the real line give rise to a vanishing $1$-form.
Write the corresponding parallel transport as
$R=(g(x,y), 1) \in TAut_2$. Integrating the equation
$\d \,\, \ch_\xi(x,y)=\omega_2(\ch_\xi(x,y))$ along the lower eyelid,
we obtain
$$
R(\ch(x,y))=\ch(\mathrm{Ad}_{g(x,y)}x,y)=\ch(y,x),
$$
and this equation implies $g(x,y)=\exp(y)$, as required.
\fin

Note that the path along the upper eyelid (oriented in the counter-clockwise
direction) can be obtained by applying the involution $\sigma_1$ to the lower
eyelid. Hence, the corresponding parallel transport is given by
$\tau_1(R)=R^{2,1}$.

\begin{prop}
The element $F \in TAut_2$ verifies the following identities,
$$
F= e^{t/2}\tau_1(F)\tau_1(R^{-1}) =e^{-t/2}\tau_1(F)R.
$$
\end{prop}

\begin{figure}[h!]
\begin{center}
\includegraphics[width=11cm]{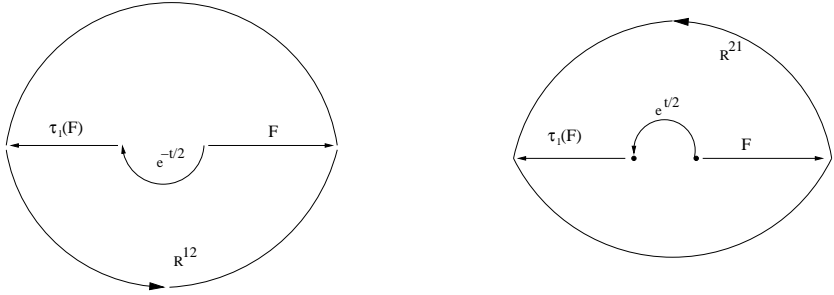}
\caption{\footnotesize The $R$ matrix}\label{Rmatrix.eps}
\end{center}
\end{figure}

\noindent\textit{Proof : }
These equations express the flatness condition for two contractible
paths shown on  Fig.\ref{Rmatrix.eps}.
\fin

These equations can be re-interpreted as the property of the parallel
transport under the involution $\tau_1$,
$$
F^{2,1}=\tau_1(F)=e^{-t/2}F R^{2,1}=e^{t/2}FR^{-1}.
$$

\section{Connection $\omega_3$ and associators}\label{PhiAT}

\subsection{Connection  $\omega_3^\infty$}
An important element of our construction is the connection
$\omega_3^\infty$ which is built using the Kontsevich technique
applied to the complex plane equipped with the Euclidean angle form.

\begin{figure}[h!]
\begin{center}
\hspace{-2cm}\includegraphics[width=6cm]{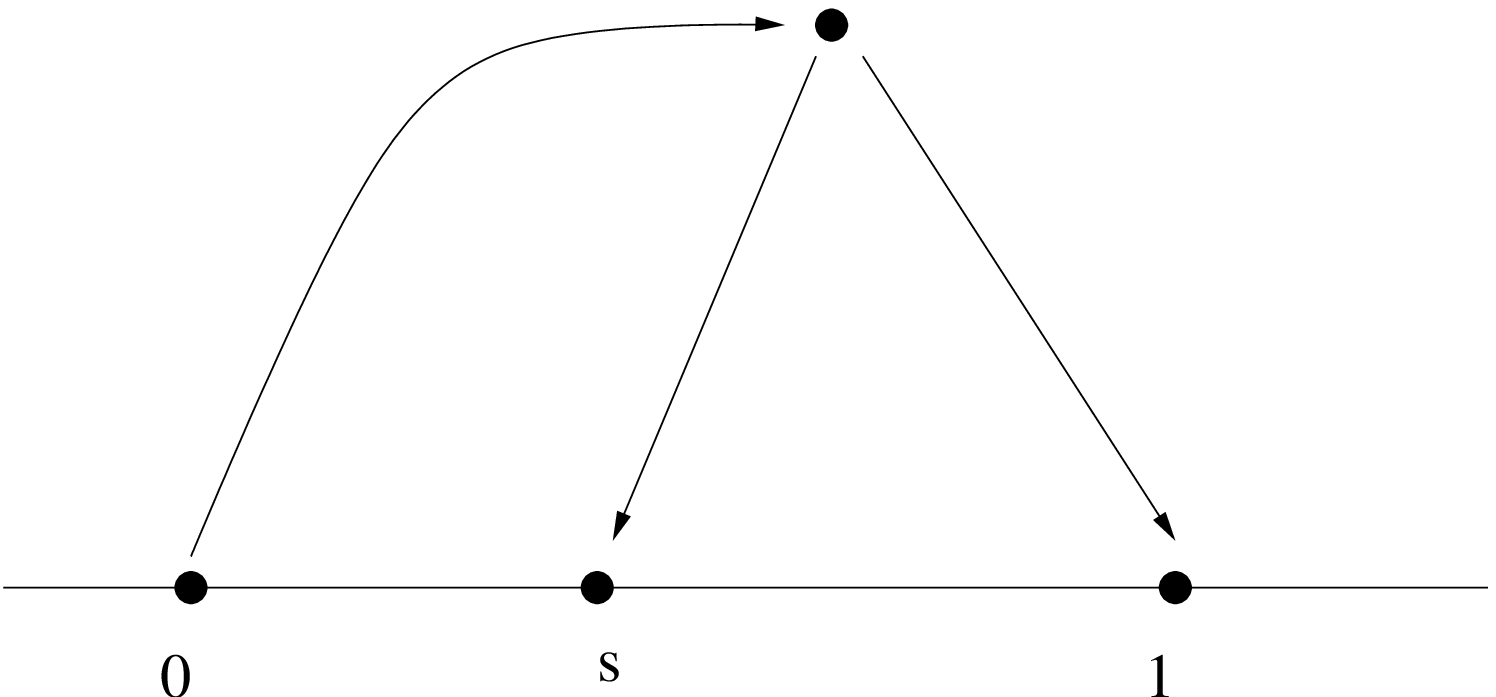}
\caption{}\label{connectioninfini.eps}
\end{center}
\end{figure}

\begin{ex}\label{exemple} Consider $\omega_3^\infty$ for 3 points
situated on the real line at the positions $0, s, 1$ (see
Fig.~\ref{connectioninfini.eps}). The Euclidean angle form is
$\d \theta$ with $\tan(\theta)= \frac{y_2-y_1}{x_2-x_1}$.
One of the simplest trees is shown on Fig.~\ref{connectioninfini.eps}.
The corresponding 3-form is given by the following expression,
\begin{multline}\frac1{1+ (\frac{y}{x})^2} \d\left(\frac{y}x\right)\wedge \frac{1}{1+ (\frac{y}{x-s})^2} \d\left(\frac{y}{(x-s)}\right) \wedge \frac1{1+(\frac{y}{x-1})^2}\d\left(\frac{y}{(x-1)}\right)=\\ - \frac{y^2}{(x^2+y^2)((x-s)^2 + y^2)( (x-1)^2+y^2)} dx\wedge dy\wedge ds.\end{multline}
By \cite{AMM} \S I.1,  the orientation is given by $- dx\wedge dy\wedge ds$,
and  one gets
$$
\omega^\infty_{\G^{(1)}}= -\frac1{8\pi^2}
\left(\frac{\log(1-s)}s  + \frac{\log(s)}{(1-s)}\right)\, \d s.
$$
This 1-form is integrable (semi-algebraic), and one has
$\int_0^1 \omega^\infty_{\G^{(1)}}= \frac1{24}.$
\end{ex}

\begin{rem}\label{remsym2} Let $\alpha: z \mapsto \bar{z}$ be a the complex
conjugation, and let $\kappa$ be an involution of $\tder_3$ defined by
formula
$$
(a(x,y,z),b(x,y,z),c(x,y,z)) \mapsto
(a(-x,-y,-z),b(-x,-y,-z),c(-x,-y,-z)).
$$
Then, $\alpha^* \omega_3^\infty= \kappa(\omega_3^\infty)$. The proof
is similar to the one of Proposition \ref{symmetries}.
\end{rem}

\begin{prop}
Connection $\omega_3^\infty$ is flat and takes values in $\krv_3$.
\end{prop}

\noindent\textit{ Proof :}
The flatness condition $\d \omega_3^\infty +
\frac 12[\omega_3^\infty, \omega_3^\infty]=0$ is obtained by replacing
the hyperbolic angle form on the upper half-plane by the Euclidean angle form
on the complex plane in the proof of Theorem \ref{theocourbure}.

In the same way, equations (\ref{KVcharles1}) and (\ref{KVcharles2})
are valid for the connection $\omega_n^\infty$.
By Kontsevich Vanishing  Lemma (Lemma 6.6 \cite{Kont}), $w_\G^\infty(\xi)=0$
for Kontsevich coefficient associated to nontrivial tree graphs.
Hence, $\ch^\infty_\xi(x,y,z)=x+y+z$.
Furthermore, Shoikhet's result \cite{Shoi} implies $\duf_\xi^\infty(x,y,z)=0$.
Therefore, differential equations (\ref{KVcharles1}) and (\ref{KVcharles2}) yield
$$
\omega_3^\infty(x+y+z)=0 \quad , \quad \div(\omega_3^\infty)=0.
$$
That is, $\omega_3^\infty$ takes values in $\krv_3$ as required.
\fin

\subsection{Associator}
We will use the following notation. Recall that $T=(u,v)\in TAut_2$ is an automorphism of
$\lie_2$ acting by
$$
(x, y) \mapsto \left(\mathrm{Ad}_u  x, \mathrm{Ad}_v y \right) .
$$
We denote $T^{1,2}=(u(x,y),v(x,y),1) \in TAut_3$,
$T^{12,3}=(u(x+y,z),u(x+y,z),v(x+y,z) \in TAut_3$ {\em etc.}
For $F \in TAut_2$ the parallel transport from the iris to the right
corner of the eye $\ov{C}_{2,0}$, we define
$$
\Phi = F^{1, 23} F^{23}(F^{12,3}F^{12})^{-1}  \in TAut_3.
$$
This element is the main topic of study in this Section.

\begin{prop}\label{proptransport} The element $\Phi$ coincides with the parallel transport
for the equation  $\d g= - g \omega_3^\infty$  between positions $1(23)$ to $(12)3$.
\end{prop}

\noindent \textit{Proof : }  Consider the following path in the configuration space
$\ov{C}_{3,0}$:

First, place $z_1,z_2,z_3$ on the stratum at infinity and move them
along the horizontal line (the real axis of the complex plane at infinity)
from the position  $1(23)$ ($z_2$ and $z_3$ collapsed) to the position $(12)3$
($z_1$ and $z_2$ collapsed). The connection at infinity is $\omega_3^\infty$, and
we denote the corresponding parallel transport by $\Phi^\infty$.

Next, make $z_3$ descend from the stratum at infinity to plus infinity of the real
axis of the upper half-plane (this corresponds to moving to the right corner
of the eye for the points (12) and 3). On this stratum, the connection is $\omega_2^{12,3}$,
and the parallel transport is given by $F^{12,3}$.

Continue with descending both $z_1$ and $z_2$ to the real axis. The connection
on this stratum is $\omega_2^{1,2}$, and the parallel transport gives $F^{1,2}$.

Then, move $z_2$ to the vicinity of $z_3$ along the real axis of the upper half-plane.
The parallel transport is trivial since the connection vanishes along the real axis.

Finally,  lift $z_2$ and $z_3$ from the real axis and make them collapse on each other
(parallel transport $(F^{2,3})^{-1}$), and lift $z_1$ from the real axis
and make it collapse with  $z_2=z_3$ (parallel transport $(F^{1,23})^{-1}$).

Thus, we made a loop and returned to the position $1(23)$ at infinity. This loop
is contractible, and the total parallel transport is trivial
by the flatness property of the connection.
Hence,
$$
\Phi^\infty  F^{12,3} F^{1,2}(F^{2,3})^{-1}(F^{1,23})^{-1} = 1,
$$
and we obtain
$$
\Phi^\infty=F^{1,23} F^{2,3} (F^{12,3} F^{1,2})^{-1} = \Phi,
$$
as required.
\fin

We have  $\Phi= F^{1,23}F^{23} (F^{12, 3}F^{12})^{-1}$,
and the first term of $\Phi$ is given by
$$
\Phi(x,y,z)= 1 - \frac1{24} ([y,z], -[x,z],[y,z]) +\ldots=
1 + \frac1{24}[t^{1,2}, t^{2,3}] + \dots ,
$$
with $t^{1,2}=(y,x,0)$ and $t^{2,3}=(0, z,y)$. Here we used the fact that
$\int_0^1 \omega^\infty_{\G^{(1)}}=\frac1{24}$ (see Example considered above),
and the minus sign is coming from the orientation of the boundary stratum
$C_{3} \subset \partial{C_{3, 0}}$. The main properties of the element $\Phi$
are summarized in the following Theorem.

\begin{theo}\label{theoassociator}
The element $\Phi$ satisfies associator axioms.
\end{theo}

\noindent \textit{Proof : }
The axioms to verify are as follows: $\Phi$ is a group like element and it is a solution
of the following equations,

\begin{eqnarray}\nonumber
\Phi^{3,2,1}\Phi^{1,2,3}=1 \quad (i) \\\nonumber
\Phi^{1,2, 34} \Phi^{12, 3, 4}= \Phi^{2, 3, 4} \Phi^{1, 23, 4} \Phi^{1, 2, 3}\quad (ii) \\\nonumber
\exp\left(\pm \frac12 t_{12}\right) \Phi^{3,1,2} \exp\left(\pm \frac12 t_{13}\right) \Phi^{2, 3, 1} \exp\left(\pm\frac12 t_{23} \right) \Phi^{1,2,3} =\\\exp\left(\pm\frac12 (t_{12}+t_{13}+ t_{23})\right)\quad (iii)\nonumber
\end{eqnarray}

(i) - $(\Phi^{1,2,3})^{-1}$ is the parallel transport between positions $1(23)$ to $(12)3$.
Let $\beta$ be the reflection with respect to the vertical axis. Then, similar to
Proposition  \ref{symmetries}, we obtain $\beta^* \omega_3^\infty = (\omega_3^\infty)^{3,2,1}$.
Hence, we get $ \Phi^{3,2,1}= (\Phi^{1,2,3})^{-1}$, as required.\\

(iii) - Consider  the following path:
$$
1(23) \mapsto (12)3 \mapsto (21)3 \mapsto 2(13) \mapsto 2(31)  \mapsto (23) 1 \mapsto 1(23).
$$
Here the last step is by moving the collapsed pair $(23)$ around the point $1$ along the iris of the
corresponding $\ov{C}_{2,0}$ stratum. By the flatness property, the total
parallel transport is trivial. This gives exactly the pair of hexagonal equations \textit{(iii)}
by using the equation (i) and the fact that $c=t_{12}+t_{13}+ t_{23}$ is central in
$\sder_3$ (see Proposition 3.4, \cite{AT1}). The plus or minus sign
in equation \textit{(iii)} depends
on the choice of the (clockwise or anti-clockwise) semi-circle
for each exchange of two points ({\em e.g.} 1 moves above 2 or below 2 in
the move $(12)3 \mapsto (21)3$).

\begin{figure}[h!]
\begin{center}
\includegraphics[width=8cm]{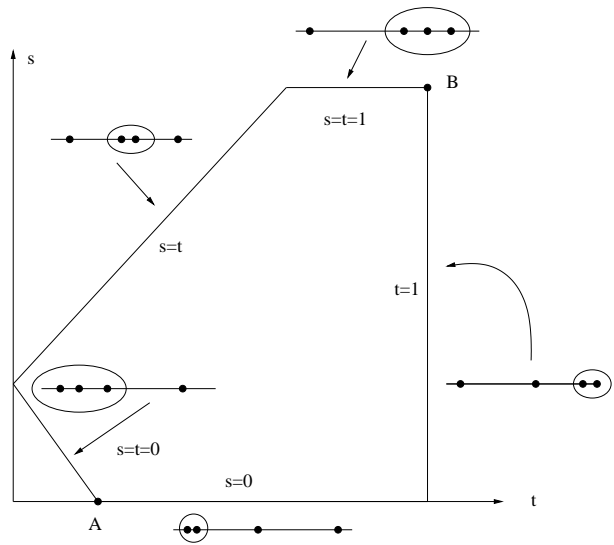}
\caption{\footnotesize The compactification space of $4$ positions on a line}\label{4position.eps}
\end{center}
\end{figure}

(ii) -  Consider four points $z_1=0, z_2=s, z_3=t, z_4=1$ on the horizontal line
(the real axis of the complex plane) representing
a point of the configuration space of the complex plane placed at infinity of $\ov{C}_{4,0}$.
The path
$$
((12)3)4 \mapsto(1(23))4  \mapsto 1((23)4) \mapsto 1(2(34))
$$
is contractible. Hence, the parallel transport defined by the flat connection
$\omega^\infty_4$ is trivial.
It is easy to see that it reproduces the pentagon equation \textit{(ii)} (see Fig.~ \ref{4position.eps}).

\fin

Note that $\Phi$ is an element of the group $KRV_3$ which contains the subgroup $T_3=\exp(\t_3)$.
If $\Phi$ is actually an element of $T_3$, it becomes a Drinfeld associator. Since $\Phi$
is even ($\kappa(\Phi)=\Phi$, by Remark~\ref{remsym2}), it does not coincide the
Knizhnik-Zamolodchikov associator (the only
known Drinfeld associator defined by an explicit formula).

\begin{conj} The element $\Phi$ is a Drinfeld associator.
That is, $\Phi \in T_3 \subset KRV_3$.
\end{conj}

In a recent work \cite{SW}, Severa and Willwacher prove this conjecture
affirming that the element $\Phi$ is indeed a new Drinfeld associator admitting
a presentation as a parallel transport of the flat connection $\omega_3^\infty$ defined by
explicit formulas.


\begin{thebibliography}{99}
\small


\bibitem{AM2} Alekseev, A.; Meinrenken, E.,  \textit{ On the Kashiwara-Vergne conjecture.}
 Invent. Math.  \textbf{164}  (2006), no. 3, 615--634.
\bibitem{AMM}
Arnal, D.; Manchon, D. ; Masmoudi, M., Choix des signes pour la
formalit\'e de {K}ontsevich.
 \textit{ Pacific J. Math.} \textbf{203} (2002), 23--66.

\bibitem{AT1} Alekseev, A.; Torossian, C., \textit{ The Kashiwara-Vergne conjecture and Drinfeld's associators.} arxiv: 0802.4300

\bibitem{AST} Andler, M.; Sahi, S.; Torossian, C., Convolution of invariant distributions:
 proof of the Kashiwara-Vergne conjecture. \textit{ Lett. Math. Phys.} \textbf{ 69}  (2004), 177--203.

\bibitem{CKT} Cattaneo, A.S.;  Keller, B.;
Torossian, C.;  Brugui\`eres, A., \textit{D\'eformation,
quantification, th\'eorie de Lie.} Collection Panoramas et Synth\`ese
no. \textbf{20}, SMF, 2005.

\bibitem{Dr} Drinfeld, V. G., \textit{On quasitriangular quasi-Hopf algebras and on a group that is closely connected with $\mathrm{Gal}(\overline{Q}/Q)$}.
(Russian)  Algebra i Analiz  2  (1990),  no. 4, 149--181;  translation in  Leningrad Math. J.  2  (1991),  no. 4, 829--860
\bibitem{FW} Felder, G. ; T. Willwacher,  \textit{On the (ir)rationality of Kontsevich weights} arXiv: 0808.2762.

\bibitem{KV} Kashiwara, M. ;  Vergne, M., The Campbell-Hausdorff
formula and invariant hyperfunctions. \textit{Inventiones Math.
}\textbf{47} (1978), 249--272.

\bibitem{Ka} Kathotia, V., \textit{Kontsevich's universal formula for deformation
  quantization and the {C}ampbell-{B}aker-{H}ausdorff formula}. Internat. J.
  Math. \textbf{11} (2000), no.~4, 523--551.
\bibitem{Kont} Kontsevich, M.,   \emph{Deformation quantization of {P}oisson
manifolds, {I}}, Preprint of the IH{\'E}S, October 1997,
q-alg/9709040 published in   \textit{Lett. Math.Phys.} \textbf{66
}(2003), no. 3, 157--216.
\bibitem{SW} Severa, P. ; Willwacher, T, Equivalence of formalities of the little discs operad,
preprint  arXiv:0905.1789 .
\bibitem{Shoi} Shoikhet, B.,  Vanishing of the Kontsevich integrals of the wheels. EuroConf\'erence Moshe Flato 2000, Part II (Dijon).  Lett. Math. Phys.  \textbf{56 } (2001),  no. 2, 141--149.
\bibitem{To1}  Torossian, C.,  Sur la conjecture
combinatoire de Kashiwara-Vergne.  \textit{J. Lie Theory}
\textbf{12} (2002),  no. 2, 597--616.

\end{thebibliography}
\end{document}